\documentclass[11pt,english]{article}
\usepackage[T1]{fontenc}
\usepackage[utf8]{inputenc}
\usepackage{lmodern}
\usepackage{microtype}
\usepackage{titling}
\usepackage{amssymb,latexsym,amsmath,amsthm,amsfonts,mathtools}
\usepackage{mathrsfs}
\usepackage{graphicx}
\usepackage{relsize}
\usepackage{fancyhdr}
\usepackage{pgfplots}
\usepackage{tikz,tikz-cd}
\usetikzlibrary{arrows.meta,positioning,shapes.geometric,shapes,backgrounds}
\usepackage{enumerate,enumitem}
\usepackage{authblk}
\usepackage[a4paper,margin=26mm]{geometry}
\usepackage{framed}
\usepackage{booktabs}
\usepackage[colorlinks=true,linkcolor=black,citecolor=black,urlcolor=black]{hyperref}
\allowdisplaybreaks[2]
\newcommand{\g}{\mathfrak{g}}
\newcommand{\ten}{\otimes}
\newcommand{\id}{\mathrm{id}}
\newcommand{\eps}{\varepsilon}
\newcommand{\Del}{\Delta}
\newcommand{\act}{\mathbin{\diamond}}
\newcommand{\tri}{\mathbin{\triangleright}}
\newtheorem{lemma}{Lemma}
\newtheorem{proposition}[lemma]{Proposition}
\newtheorem{theorem}[lemma]{Theorem}

\newtheorem{definition}[lemma]{Definition}

\theoremstyle{definition}
\newtheorem{remark}[lemma]{Remark}
\newcounter{wordlabel}

\tikzset{box/.style={draw,thick,rounded corners,align=center,text width=5cm,inner sep=6pt},arrow/.style={-{Stealth[length=2mm]},thick},alabel/.style={font=\small\itshape,fill=white,inner sep=1pt}}

\begin{document}
\title{\textbf{Post-Hopf structures on $U_q(\mathfrak{sl}_2)$ and post-Lie bialgebras}}
\author{Andrea Rivezzi}
\affil{{\sl Department of Mathematics, Faculty of Science, University of Zagreb}.\\
{\sl Bijeni\v{c}ka 30, 10000 Zagreb, Croatia}}
\date{}
\maketitle
\vspace{-1cm}
\begin{abstract}
\noindent
We classify post-Hopf structures on the Drinfeld–Jimbo quantum group \(U_q(\mathfrak{sl}_2)\), for \(q\in\mathbb C^\times\) not a root of unity, proving that there are exactly three isomorphism classes. We then investigate the semiclassical limits of post-Hopf quantized universal enveloping algebras and derive a compatibility condition between the induced post-Lie product and Lie cobracket. This leads to the notion of a post-Lie bialgebra. We construct a nontrivial topological post-Hopf structure on \(U_\hbar(\mathfrak{sl}_2)\) and determine its semiclassical limit. For the standard Lie bialgebra structure on \(\mathfrak{sl}_2\), we classify all compatible post-Lie products up to isomorphism and show that the nontrivial class admits a quantization. Finally, we formulate the general quantization problem for post-Lie bialgebras.
\end{abstract}
\tableofcontents
\section*{Introduction}

Quantum groups originated in the study of quantum integrable
systems and have since become fundamental objects in several
areas of mathematics and mathematical physics. One of the
first examples appearing in the literature was the quantum deformation of the universal
enveloping algebra $U(\mathfrak{sl}_2)$ introduced by Kulish
and Reshetikhin \cite{KulishReshetikhin}. This construction
was subsequently generalized by Drinfeld and Jimbo to quantum
enveloping algebras $U_q(\mathfrak g)$ associated with
semisimple Lie algebras $\mathfrak g$ \cite{Drinfeld,Jimbo}.
The passage to the semiclassical limit revealed the underlying
structure of a Lie bialgebra, namely a Lie algebra equipped with
a compatible Lie cobracket, and naturally led to the converse
question of whether every Lie bialgebra admits a quantization \cite{DrinfeldProblems}.
This problem was solved by Etingof and Kazhdan, who constructed a quantization functor from Lie bialgebras to quantized universal enveloping algebras (QUEAs) \cite{EtingofKazhdan}. 
Moreover, in the sixth paper of their series, they proved
that their quantization of a Kac--Moody Lie bialgebra with
its standard quasitriangular structure is isomorphic to
the corresponding Drinfeld--Jimbo quantum group
\cite{EtingofKazhdanVI}.

Over the past two decades, considerable attention has also
been devoted to post-Lie algebras, which are Lie algebras
equipped with an additional bilinear product $\tri$ satisfying compatibility identities that, in particular,
make $[x,y]_{\tri}=[x,y]+x\tri y-y\tri x$
another Lie bracket.
Introduced by Vallette in an operadic context
\cite{Vallette}, these structures have attracted interest
from several mathematical communities through their
connections with differential geometry and geometric
numerical integration \cite{MuntheKaasLundervold},
Rota--Baxter operators and the Yang--Baxter equation
\cite{BaiGuoNi}, and regularity structures arising in
the theory of stochastic partial differential equations
\cite{BrunedKatsitsiadis}.
More recently, Li, Sheng, and Tang introduced post-Hopf
algebras \cite{LST}, providing a Hopf-algebraic counterpart
to post-Lie algebras. In particular, every post-Lie algebra
induces a post-Hopf structure on its universal enveloping
algebra, while the primitive elements of a post-Hopf algebra
inherit a post-Lie product. These constructions establish
an adjunction extending the classical universal enveloping
algebra--primitive elements correspondence.

In view of these developments, it is natural to ask whether
the existing quantization functors for Lie bialgebras can be
extended to quantize simultaneously a Lie cobracket and a
post-Lie product on the same Lie algebra. This requires,
first of all, understanding which compatibility conditions
these two structures must satisfy in order to arise together
as the semiclassical limit of a post-Hopf QUEA.
A first step in this direction was taken by the author,
together with Sciandra and Weber, in \cite{RSW}, where
infinitesimal deformations of post-Lie and post-Hopf algebras
and their compatibility with the universal enveloping
algebra--primitive elements correspondence were investigated.

The purpose of the present paper is to lay the groundwork
for a theory of quantization of post-Lie bialgebras,
starting from explicit quantum-group constructions and
developing the appropriate semiclassical framework.
We begin by classifying all post-Hopf structures on the
Drinfeld--Jimbo quantum group $U_q(\mathfrak{sl}_2)$,
for $q\in\mathbb C^\times$ not a root of unity.
We prove that there are exactly three isomorphism classes,
represented by $\act_+$, $\act_-$, and $\act_q$, see Theorem \ref{thm:classification}.
Their values on generators are summarized in the
following tables:
\[
\begin{array}{c|cccc}
\act_+ & K & K^{-1} & E & F\\ \hline
K      & K & K^{-1} & E & F\\
K^{-1} & K & K^{-1} & E & F\\
E      & 0 & 0 & 0 & 0\\
F      & 0 & 0 & 0 & 0
\end{array}
\qquad
\begin{array}{c|cccc}
\act_- & K & K^{-1} & E & F\\ \hline
K      & K & K^{-1} & -E & -F\\
K^{-1} & K & K^{-1} & -E & -F\\
E      & 0 & 0 & 0 & 0\\
F      & 0 & 0 & 0 & 0
\end{array}
\qquad
\begin{array}{c|cccc}
\act_q & K & K^{-1} & E & F\\ \hline
K      & K & K^{-1} & q^{-2}E & q^2F\\
K^{-1} & K & K^{-1} & q^2E & q^{-2}F\\
E      & 0 & 0 & E & -q^2F\\
F      & 0 & 0 & -qE & q^{-1}F
\end{array}
\]

We next study the semiclassical limit of a post-Hopf QUEA. We prove that the induced post-Lie product $\tri$ and Lie cobracket $\delta$ on the underlying Lie algebra $\mathfrak g$ satisfy the compatibility condition

\[
\delta(x\tri y)
=
(D_x\otimes\id+\id\otimes D_x)\delta(y),
\qquad D_x(y)=x\tri y,
\]

for all $x,y\in\mathfrak g$. This motivates our definition of a \emph{post-Lie bialgebra}: a Lie bialgebra equipped with a post-Lie product satisfying the above identity, see Theorem \ref{thm:semiclassical-compatibility} and Definition \ref{def:postlie-bialgebra}. 

We then construct a nontrivial topological post-Hopf structure on the formal quantum group $U_\hbar(\mathfrak{sl}_2)$ and determine its semiclassical limit. Writing $e,f,h$ for the standard generators of $\mathfrak{sl}_2$, the resulting post-Lie bialgebra structure is given by

\[
\begin{gathered}
\delta(h)=0,\qquad
\delta(e)=h\wedge e,\qquad
\delta(f)=h\wedge f,\\[3pt]
x\tri y=[\varphi(x),y],
\qquad
\varphi(e)=\frac12h,\quad
\varphi(f)=-\frac12h,\quad
\varphi(h)=-h.
\end{gathered}
\]
We also classify all post-Lie products on
$\mathfrak{sl}_2$ compatible with the standard
Lie cobracket (Theorem~\ref{thm:sl2-plb}),
complementing the classification of post-Lie
structures in \cite{PanLiuBaiGuo}.
Up to Lie bialgebra automorphisms, there are
exactly two isomorphism classes: the trivial one
and the nontrivial one displayed above.
Both classes admit quantizations, the nontrivial
one being realized by our construction on
$U_\hbar(\mathfrak{sl}_2)$.

In this sense, our approach parallels the early development of quantum groups, where the deformation of $U(\mathfrak{sl}_2)$ provided a concrete starting point for the theory. Motivated by these examples and by the classical quantization problem for Lie bialgebras, we conclude by formulating the general problem of quantizing post-Lie bialgebras.

\medskip

\noindent\textit{Structure of the paper.} 
Section 1 recalls the structure of $U_q(\mathfrak{sl}_2)$, constructs explicit post-Hopf structures, and proves their classification. Section 2 introduces post-Lie bialgebras through semiclassical limits, constructs a nontrivial post-Hopf structure on $U_\hbar(\mathfrak{sl}_2)$, classifies compatible post-Lie products on $\mathfrak{sl}_2$, and discusses the general quantization problem.

\medskip

\noindent\textit{Acknowledgments.} 
The author is supported by the Croatian Science Foundation (HRZZ), 
grant no. IP-2025-02-6514 and is a member of the ``National Group for Algebraic and Geometric Structures and their Applications'' (GNSAGA-INdAM). 

\medskip

\noindent\textit{AI disclosure.} 
The author used OpenAI's ChatGPT (GPT-6) for
language editing, proofreading, and assistance
in reviewing the exposition and consistency of
mathematical arguments. All AI-assisted suggestions
were independently verified by the author, who
takes full responsibility for the mathematical
results and final content of the manuscript.

\section{Post-Hopf structures on $U_q(\mathfrak{sl}_2)$}
Throughout this section, all vector spaces are over $\mathbb{C}$, and $\ten$ denotes the tensor product of complex vector spaces. 
We use Sweedler notation
$\Del(a)=a_1\otimes a_2$,
with summation understood, and adopt the
comultiplication convention of \cite{Jantzen}.
\subsection{Definition and graded structure of $U_q(\mathfrak{sl}_2)$}\label{subsec:uq}
Let $q\in\mathbb{C}^\times$ not be a root of unity.  Recall that the Drinfeld--Jimbo quantum group $U_q(\mathfrak{sl}_2)$
is the algebra generated by elements
$E,F,K,K^{-1}$ subject to the relations
\begin{equation}\label{eq:uq-rel}
KK^{-1}=K^{-1}K=1,\quad KEK^{-1}=q^2E,\quad KFK^{-1}=q^{-2}F,\quad
EF-FE=\frac{K-K^{-1}}{q-q^{-1}}.
\end{equation}
Its Hopf algebra structure is given by the comultiplication
\begin{equation}\label{eq:uq-cop}
\Del(K)=K\otimes K,\qquad
\Del(E)=E\otimes 1+K\otimes E,\qquad
\Del(F)=F\otimes K^{-1}+1\otimes F.
\end{equation}
The counit is $\varepsilon(K)=1$, $\varepsilon(E)= \varepsilon(F)=0$, and the antipode is $S(K)=K^{-1}$, $S(E)=-K^{-1}E$, and $S(F)=-FK$.

It is well known that $U_q(\mathfrak{sl}_2)$ admits the  Poincar\'e--Birkhoff--Witt basis $\{ F^rK^kE^s\}_{r,s\ge0,\ k\in\mathbb{Z}}$.
Also, relations \eqref{eq:uq-rel} are homogeneous with respect to the
assignment
\[
\deg(E)=1,\qquad \deg(F)=-1,\qquad
\deg(K^{\pm1})=0.
\]
Consequently, they induce a $\mathbb{Z}$-grading
\begin{equation}\label{eq:grading}
U_q(\mathfrak{sl}_2)
=\bigoplus_{n\in\mathbb{Z}} U_q(\mathfrak{sl}_2)_n,
\end{equation}
where $U_q(\mathfrak{sl}_2)_n$ is spanned by the
PBW monomials $F^rK^kE^s$ such that $s-r=n$.
In particular, for every homogeneous element
$a\in U_q(\mathfrak{sl}_2)_n$, one has $KaK^{-1}=q^{2n}a$.
Moreover, the comultiplication is compatible with
this grading, namely
\begin{equation}\label{eq:graded-cop}
\Del\bigl(U_q(\mathfrak{sl}_2)_n\bigr)
\subseteq
\bigoplus_{i+j=n}
U_q(\mathfrak{sl}_2)_i
\otimes U_q(\mathfrak{sl}_2)_j.
\end{equation}

We shall use the following classification of Hopf
algebra automorphisms of $U_q(\mathfrak{sl}_2)$;
see \cite[\S 6.1.6]{KlimykSchmudgen}.

\begin{proposition}
\label{prop:aut}
Every Hopf algebra automorphism of
$U_q(\mathfrak{sl}_2)$ is uniquely of the form
\begin{equation}\label{eq:phit}
\Phi_t(K)=K,\qquad
\Phi_t(E)=tE,\qquad
\Phi_t(F)=t^{-1}F,
\qquad t\in\mathbb{C}^\times.
\end{equation}
\end{proposition}
Next, for grouplike elements $g,h\in U_q(\mathfrak{sl}_2)$,
let $P_{g,h}
=\{x\in U_q(\mathfrak{sl}_2):
\Del(x)=x\otimes h+g\otimes x\}$
denote the space of $(g,h)$-skew-primitive elements.
\begin{lemma}
We have
\begin{align}
P_{K,K}&=0,
\label{eq:skew-KK}\\
P_{K,1}
&=\operatorname{span}_\mathbb{C}\{E,FK,K-1\},
\label{eq:skew-K1}\\
P_{1,K^{-1}}
&=\operatorname{span}_\mathbb{C}\{F,EK^{-1},1-K^{-1}\}.
\label{eq:skew-1Kinv}
\end{align}
\end{lemma}
\begin{proof}
\eqref{eq:skew-K1} follows from
\cite[\S 3.1.2, proof of Proposition~6]{KlimykSchmudgen}\footnote{The different
coproduct convention used in the cited reference has been taken into account.}. Moreover, right multiplication
by $K^{-1}$ gives
$P_{1,K^{-1}}=P_{K,1}K^{-1}$, proving
\eqref{eq:skew-1Kinv}.
To show \eqref{eq:skew-KK}, note that if $r+s>0$, the coproduct of $F^rK^kE^s$ contains
the terms $F^rK^kE^s\otimes K^{k-r}$ and $K^{k+s}\otimes F^rK^kE^s$,
while every other term contains an $E$ or $F$ in
both tensor factors.
Thus, by linear independence of tensor products
of PBW monomials, if $a=\sum c_{r,k,s}F^rK^kE^s$
is primitive, then $c_{r,k,s}\neq0$ and $r+s>0$
imply $k-r=k+s=0$, which is a contradiction.
Hence $a\in\mathbb{C}[K^{\pm1}]$, which has no nonzero
primitive elements since $\Del(K^k)=K^k\otimes K^k$.
Finally, since $P_{K,K}=P_{1,1}K$ and
$P_{1,1}=0$, we conclude that $P_{K,K}=0$.
\end{proof}

\subsection{Post-Hopf algebras}\label{subsec:posthopf}

Recall that a \emph{post-Hopf
algebra} is a pair $(A,\act)$, where $A$ is a Hopf algebra and $\act$ is a
coalgebra morphism $\act:A\otimes A\to A$ such that, for every $a,b,c\in A$,
\begin{align}
a\act(bc)&=(a_1\act b)(a_2\act c),
\label{eq:PH1}\\
a\act(b\act c)&=\bigl(a_1(a_2\act b)\bigr)\act c,
\label{eq:PH2}
\end{align}
and such that $\alpha_{\act}:A\to\operatorname{End}(A)$,
$\alpha_{\act}(a)(b)=a\act b$, is convolution invertible in $\mathrm{Hom}(A, \mathrm{End}(A))$. 
The convolution condition means that there exists a unique
$\beta_{\act}:A\to\operatorname{End}(A)$ such that
\begin{equation}\label{eq:beta}
\alpha_{\act}(a_1)\circ\beta_{\act}(a_2)
=\beta_{\act}(a_1)\circ\alpha_{\act}(a_2)
=\eps(a)\id_A.
\end{equation}
From the post-Hopf axioms, it follows that 
\begin{equation}\label{eq:norm}
1\act a=a,\qquad a\act1=\eps(a)1,\qquad
S(a\act b)=a\act S(b).
\end{equation}
Every Hopf algebra admits the trivial post-Hopf structure  $a\act_+ b=\eps(a)b$. In the case of $U_q(\mathfrak{sl}_2)$, its values on generators are 
\[
\begin{array}{c|cccc}
\act_+ & K & K^{-1} & E & F\\ \hline
K      & K & K^{-1} & E & F\\
K^{-1} & K & K^{-1} & E & F\\
E      & 0 & 0 & 0 & 0\\
F      & 0 & 0 & 0 & 0
\end{array}
\]

As we shall see in the proof of Theorem~\ref{thm:classification},
every post-Hopf structure on $U_q(\mathfrak{sl}_2)$
is uniquely determined by its values on pairs
of generators.

Two post-Hopf structures $\act$ and $\act'$ on the same
Hopf algebra $A$ are said to be isomorphic if there exists
a Hopf algebra automorphism $\Psi$ of $A$ such that $\Psi(a\act b)=\Psi(a)\act'\Psi(b)$
for all $a,b\in A$.

\subsubsection{Post-Hopf structures on $\mathbb{Z}$-graded Hopf algebras}

We now give a construction of post-Hopf
structures on graded Hopf algebras, which will be useful in the
sequel.

\begin{proposition}\label{prop:graded}
Let $A=\bigoplus_{n\in\mathbb{Z}}A_n$ be a $\mathbb{Z}$-graded Hopf algebra
and let $\{\chi_n:A\to\mathbb{C}\}_{n\in\mathbb{Z}}$ be a 
family of linear maps satisfying, for all
$m,n,r\in\mathbb Z$, $a\in A$, and $b\in A_r$,
the following conditions:
\begin{equation}\label{eq:chi-axioms}
\begin{gathered}
\chi_0=\eps,\qquad \chi_n(1)=1,\qquad
\chi_m*\chi_n=\chi_{m+n}, \qquad
\chi_m(ab)=\chi_{m-r}(a)\chi_m(b)
\end{gathered}
\end{equation}
where $*$ denotes the convolution product on
$\operatorname{Hom}(A,\mathbb{C})$. For $a\in A$ and $b\in A_n$, define
$a\act b=\chi_n(a)b$ and extend by bilinearity.
Then $\act$ is a post-Hopf structure on $A$, with
convolution-inverse action map given by $\beta_{\act}(a)(b)=\chi_{-n}(a)b$.
\end{proposition}

\begin{proof}
Define $\act$ as in the statement. We first prove that $\act$ is a coalgebra morphism. 
For homogeneous $b\in A_n$, write $\Del(b)=\sum_{r+s=n}\Del_{r,s}(b)$, where $
\Del_{r,s}(b)\in A_r\otimes A_s$.
Then, using the properties of the maps
$\chi_m$, we obtain
\begin{align*}
\Del(a\act b)
=\Del(\chi_n(a)b)
&=\chi_n(a)\Del(b)=\sum_{r+s=n}\chi_{r+s}(a)\Del_{r,s}(b)\\&=\sum_{r+s=n}\chi_r(a_1)\chi_s(a_2)
\Del_{r,s}(b)=(a_1\act b_1)\otimes(a_2\act b_2).
\end{align*}
Moreover, since $\chi_0=\eps$ and the counit of a graded Hopf algebra vanishes
on homogeneous components of nonzero degree, we have
$\eps(a\act b)=\eps(a)\eps(b)$.
Thus, $\act$ is a coalgebra morphism.
Next,  we verify \eqref{eq:PH1},\eqref{eq:PH2}. Let $b\in A_r$ and $c\in A_s$. By
\eqref{eq:chi-axioms}, we have
\[
\begin{aligned}
a\act(bc)
=\chi_{r+s}(a)bc=\chi_r(a_1)\chi_s(a_2)bc=(a_1\act b)(a_2\act c),
\end{aligned}
\]
which proves \eqref{eq:PH1}.
Similarly, using the last identity in
\eqref{eq:chi-axioms}, we obtain
\[
\begin{aligned}
a\act(b\act c)&=\chi_s(a)\chi_s(b)c=(\chi_{s-r}*\chi_r)(a)\chi_s(b)c=\chi_{s-r}(a_1)\chi_r(a_2)\chi_s(b)c\\
&=\chi_r(a_2)\chi_s(a_1b)c=\bigl(a_1(a_2\act b)\bigr)\act c,
\end{aligned}
\]
which proves \eqref{eq:PH2}.
Finally, for $b\in A_n$, we have
\[
\begin{aligned}
\alpha_{\act}(a_1)\bigl(\beta_{\act}(a_2)(b)\bigr)
=\chi_n(a_1)\chi_{-n}(a_2)b=(\chi_n*\chi_{-n})(a)b=\eps(a)b.
\end{aligned}
\]
The opposite convolution identity follows similarly,
since $\chi_{-n}*\chi_n=\chi_0=\eps$.
Hence $\beta_{\act}$ is the convolution inverse of
$\alpha_{\act}$, completing the proof.

\end{proof}

\subsection{Classification of post-Hopf structures on $U_q(\mathfrak{sl}_2)$}\label{subsec:class}

The aim of this section is to classify, up to isomorphism, all post-Hopf structures on
$U_q(\mathfrak{sl}_2)$. Besides the trivial one, we shall construct two nontrivial post-Hopf structures $\act_-, \act_q$ and
then prove that these exhaust all possibilities.

\subsubsection{The post-Hopf structure $\act_{-}$}

Let $C_2=\{1,g\}$ denote the cyclic group of order two and consider the
Hopf algebra morphism
\begin{equation}\label{eq:quotient}
\pi:U_q(\mathfrak{sl}_2)\longrightarrow\mathbb{C}[C_2],\qquad
\pi(E)=\pi(F)=0,\quad \pi(K)=\pi(K^{-1})=g.
\end{equation}

Consider the left action $\mathbb{C}[C_2]\otimes U_q(\mathfrak{sl}_2)
\to U_q(\mathfrak{sl}_2)$ defined by
$1\cdot b=b$ and $g\cdot b=\Phi_{-1}(b)$.
For all $a,b\in U_q(\mathfrak{sl}_2)$, set
\begin{equation}\label{eq:sign}
a\act_- b=\pi(a)\cdot b.
\end{equation}

\begin{theorem}
The pair $(U_q(\mathfrak{sl}_2), \act_-)$ is a post-Hopf algebra.
\end{theorem}
\begin{proof}
Since $1$ and $g$ are grouplike and act by
coalgebra automorphisms, for every
$x\in\mathbb{C}[C_2]$ and $b\in U_q(\mathfrak{sl}_2)$
we have $\Del(x\cdot b)
=(x_1\cdot b_1)\otimes(x_2\cdot b_2)$.
Thus, using the fact that $\pi$ is a coalgebra
morphism, we obtain
\[
\begin{aligned}
\Del(a\act_- b)
=\Del(\pi(a)\cdot b)&=(\pi(a)_1\cdot b_1)\otimes
  (\pi(a)_2\cdot b_2)\\
  &=(\pi(a_1)\cdot b_1)\otimes
  (\pi(a_2)\cdot b_2)=(a_1\act_- b_1)\otimes(a_2\act_- b_2).
\end{aligned}
\]
Moreover, $\eps(a\act_- b)
=\eps(\pi(a)\cdot b)
=\eps(\pi(a))\eps(b)
=\eps(a)\eps(b),$
hence $\act_-$ is a coalgebra morphism.
Next, we prove that $\act_-$ satisfies \eqref{eq:PH1},\eqref{eq:PH2}. 
Since $1$ and $g$ act by algebra
automorphisms, we have $x\cdot(bc)=(x_1\cdot b)(x_2\cdot c)$
for all $x\in\mathbb{C}[C_2]$ and
$b,c\in U_q(\mathfrak{sl}_2)$.
Taking $x=\pi(a)$ and using again that $\pi$
is a coalgebra morphism, we obtain \eqref{eq:PH1}.
Next, since $\pi\circ\Phi_{-1}=\pi$, we have
$\pi(a\act_-b)=\eps(a)\pi(b)$ for all
$a,b\in U_q(\mathfrak{sl}_2)$.
Consequently,
$\pi\bigl(a_1(a_2\act_-b)\bigr)=\pi(a)\pi(b)$,
and hence
\[
\begin{aligned}
\bigl(a_1(a_2\act_-b)\bigr)\act_-c
=(\pi(a)\pi(b))\cdot c=\pi(a)\cdot(\pi(b)\cdot c)=a\act_-(b\act_-c),
\end{aligned}
\]
which proves \eqref{eq:PH2}. Finally, one checks that the map
$\beta_{\act_-}:U_q(\mathfrak{sl}_2)
\to\operatorname{End}(U_q(\mathfrak{sl}_2))$
defined by $\beta_{\act_-}(a)(b)=\pi(S(a))\cdot b$
is the convolution inverse of $\alpha_{\act_-}$.
Therefore, \eqref{eq:sign} is a post-Hopf structure
on $U_q(\mathfrak{sl}_2)$. 
\end{proof}
A direct computation shows that the values of $\act_-$ on generators are
\begin{equation}\label{eq:sign-gen}
\begin{array}{c|cccc}
\act_- & K & K^{-1} & E & F\\ \hline
K      & K & K^{-1} & -E & -F\\
K^{-1} & K & K^{-1} & -E & -F\\
E      & 0 & 0 & 0 & 0\\
F      & 0 & 0 & 0 & 0
\end{array}
\end{equation} 
\subsubsection{The post-Hopf structure $\act_{q}$}

Next, we construct a nontrivial post-Hopf structure
on $U_q(\mathfrak{sl}_2)$ using the $\mathbb{Z}$-grading
introduced in Subsection~\ref{subsec:uq} and the construction of
Proposition~\ref{prop:graded}.
Consider the Laurent polynomial algebra
$V=\mathbb{C}[z,z^{-1}]$, with basis $v_j=z^j$,
$j\in\mathbb{Z}$. Then the formulas
\begin{equation}\label{eq:weight-rep}
\begin{aligned}
K^{\pm1}\cdot v_j&=q^{\pm2j}v_j, \qquad
E\cdot v_j&=\frac{q^{2j}-1}{q^{-2}-1}v_{j+1}, \qquad
F\cdot v_j&=\frac{1-q^{-2j}}{q-q^{-1}}v_{j-1}
\end{aligned}
\end{equation}
define a $U_q(\mathfrak{sl}_2)$-module algebra
structure on $V$ with respect to the multiplication
$v_iv_j=v_{i+j}$, i.e. one has $a \cdot (uv) = (a_1 \cdot u)(a_2 \cdot v)$ and $a \cdot 1_V = \varepsilon(a) 1_V$ for all $a \in U_q(\mathfrak{sl}_2)$ and $u,v \in V$.
Since the action of a homogeneous element of degree
$r$ sends $v_j$ to a scalar multiple of $v_{j+r}$,
we may define, for each $m\in\mathbb{Z}$, a linear functional
\begin{equation}\label{eq:chi-def}
\chi_m:U_q(\mathfrak{sl}_2)\to\mathbb{C} \qquad \text{such that} \qquad
a\cdot v_{-m}=\chi_m(a)v_{-m+r},
\qquad  \text{for all } a\in U_q(\mathfrak{sl}_2)_r,
\end{equation}
and extend it linearly to the whole algebra.

\begin{lemma}\label{lem:chi}
The family of maps $\{\chi_m\}_{m \in \mathbb{Z}}$ defined in \eqref{eq:chi-def} satisfies
\eqref{eq:chi-axioms}.
\end{lemma}

\begin{proof}
Since $v_0=1_V$,
we have $a\cdot v_0=\eps(a)v_0$ for every
$a\in U_q(\mathfrak{sl}_2)$, and then $\chi_0=\eps$.
Moreover, since $1$ acts as the identity on $V$,
we have $\chi_m(1)v_{-m}=1\cdot v_{-m}=v_{-m}$,
implying that $\chi_m(1)=1$ for every $m\in\mathbb{Z}$.
Next, by the definition of $\chi_m$ and the module
property, we have, for any $b\in U_q(\mathfrak{sl}_2)_r$, that $(ab)\cdot v_{-m}
=a\cdot(b\cdot v_{-m})
=\chi_m(b)\,a\cdot v_{-m+r}$, and therefore $\chi_m(ab)=\chi_{m-r}(a)\chi_m(b).$
Finally, the convolution identity follows from
the module-algebra property of $V$: for all
$a\in U_q(\mathfrak{sl}_2)_r$, we have
\[
\begin{aligned}
\chi_{m+n}(a)v_{-m-n+r}
=a\cdot v_{-m-n}=a\cdot(v_{-m}v_{-n})=(a_1\cdot v_{-m})(a_2\cdot v_{-n})=(\chi_m*\chi_n)(a)v_{-m-n+r},
\end{aligned}
\]
where the last equality follows from the graded
compatibility of the comultiplication.
\end{proof}
Combining Lemma~\ref{lem:chi} with
Proposition~\ref{prop:graded}, we obtain
the following post-Hopf structure on
$U_q(\mathfrak{sl}_2)$.

\begin{theorem}\label{thm:exceptional}
For $a\in U_q(\mathfrak{sl}_2)$ and homogeneous
$b\in U_q(\mathfrak{sl}_2)_n$, define $a\act_q b=\chi_n(a)b$
and extend by bilinearity. Then $\act_q$ is a
post-Hopf structure on $U_q(\mathfrak{sl}_2)$.
In particular, its values on the generators are
\begin{equation}\label{eq:table}
\begin{array}{c|cccc}
\act_q & K & K^{-1} & E & F\\ \hline
K      & K & K^{-1} & q^{-2}E & q^2F\\
K^{-1} & K & K^{-1} & q^2E & q^{-2}F\\
E      & 0 & 0 & E & -q^2F\\
F      & 0 & 0 & -qE & q^{-1}F
\end{array}
\end{equation}
\end{theorem}

\subsubsection{Exhaustiveness}

We now show that the three post-Hopf structures $\act_+, \act_-, \act_q$ defined above constitute a classification up to isomorphism.

\begin{theorem}\label{thm:classification}
Up to isomorphism, there are exactly three
post-Hopf structures on $U_q(\mathfrak{sl}_2)$,
represented by $\act_+$, $\act_-$ and $\act_q$.
\end{theorem}

\begin{proof}
Let $\act$ be a post-Hopf structure on
$U_q(\mathfrak{sl}_2)$, and write
$L_a(b)=a\act b$.
We shall determine all possible operators
$L_K,L_{K^{-1}},L_E,L_F$ and then show that
they uniquely determine $\act$.
Using that $K$ is grouplike, together with \eqref{eq:PH1}, we have
\[L_K(ab)= K \act(ab)= (K \act a)(K \act b) = L_K(a) L_K(b), \]
i.e. $L_K$ is an algebra morphism. 
Since $K$ is grouplike and $\act$ is a coalgebra morphism, $L_K$
is also a coalgebra morphism, and thus it is a Hopf algebra morphism. 
Since $K$ is grouplike and $\alpha_{\act}$
is convolution invertible, $L_K$ is invertible.
Thus, $L_K$ is a Hopf algebra automorphism, and by
Proposition \ref{prop:aut} we have that $L_K=\Phi_t$
for some $t\in\mathbb{C}^\times$.
Moreover, since $\Phi_t(K^{-1})=K^{-1}$,
equation \eqref{eq:PH2} yields $L_KL_{K^{-1}}
=L_{K(K\act K^{-1})}
=L_1=\id$, and thus $L_{K^{-1}}=\Phi_{t^{-1}}$.
Using the coalgebra structure of $U_q(\mathfrak{sl}_2)$ \eqref{eq:uq-cop}, together with
identity \eqref{eq:PH1} gives the following Leibniz-like formulas
\begin{align}
L_E(ab)&= E \act (ab) = (E \act a)(1 \act b) + (K \act a)(E \act b) =L_E(a)b+\Phi_t(a)L_E(b),
\label{eq:twisted-leibniz-E}\\
L_F(ab)&= F \act (ab)= (F \act a)(K^{-1} \act b) + (1 \act a) (F \act b)=L_F(a)\Phi_{t^{-1}}(b)+aL_F(b).
\label{eq:twisted-leibniz-F}
\end{align}
We next determine $L_E$ and $L_F$ on the
generators. 
Since $\act$ is a coalgebra morphism, using
\eqref{eq:uq-cop} and \eqref{eq:norm}, we obtain
\begin{align*}
\Del(L_E(K))&=(E\act K)\otimes(1\act K)
 +(K\act K)\otimes(E\act K)
=L_E(K)\otimes K+K\otimes L_E(K),\\
\Del(L_F(K))
&=(F\act K)\otimes(K^{-1}\act K)
 +(1\act K)\otimes(F\act K)
=L_F(K)\otimes K+K\otimes L_F(K).
\end{align*}
Thus, $L_E(K),L_F(K)\in P_{K,K}$, and hence
both vanish by \eqref{eq:skew-KK}.
Applying \eqref{eq:twisted-leibniz-E} and
\eqref{eq:twisted-leibniz-F} to $KK^{-1}=1$, we also obtain $L_E(K^{-1})=L_F(K^{-1})=0.$
Similarly, one has 
\[
\begin{aligned}
\Del(L_E(E))
&=(E\act E)\otimes(1\act1)
 +(K\act K)\otimes(E\act E)=L_E(E)\otimes1+K\otimes L_E(E),\\
\Del(L_F(E))
&=(F\act E)\otimes(K^{-1}\act1)
 +(1\act K)\otimes(F\act E)=L_F(E)\otimes1+K\otimes L_F(E),\\ 
\Del(L_E(F))
&=(E\act F)\otimes(1\act K^{-1})
 +(K\act1)\otimes(E\act F)=L_E(F)\otimes K^{-1}+1\otimes L_E(F),\\
\Del(L_F(F))
&=(F\act F)\otimes(K^{-1}\act K^{-1})
 +(1\act1)\otimes(F\act F)=L_F(F)\otimes K^{-1}+1\otimes L_F(F).
\end{aligned}
\]
Therefore, $L_E(E),L_F(E)\in P_{K,1}$ and $L_E(F),L_F(F)\in P_{1,K^{-1}}$.
Applying \eqref{eq:twisted-leibniz-E} and
\eqref{eq:twisted-leibniz-F} to
$KE=q^2EK$ and $KF=q^{-2}FK$, and using
$L_E(K)=L_F(K)=0$, we obtain
\[
\begin{aligned}
KL_E(E)&=q^2L_E(E)K, &
KL_F(E)&=q^2L_F(E)K,\\
KL_E(F)&=q^{-2}L_E(F)K, &
KL_F(F)&=q^{-2}L_F(F)K.
\end{aligned}
\] 

By \eqref{eq:skew-K1}, we may write $L_E(E)=\alpha E+\beta FK+\gamma(K-1)$
for some $\alpha,\beta,\gamma\in\mathbb{C}$.
Since $KL_E(E)K^{-1}=q^2L_E(E)$, using
the defining relations \eqref{eq:uq-rel},
we obtain
\[
\begin{aligned}
q^2\alpha E+q^{-2}\beta FK+\gamma(K-1)=q^2\alpha E+q^2\beta FK
+q^2\gamma(K-1).
\end{aligned}
\]
Since $E$, $FK$ and $K-1$ are linearly
independent and $q$ is not a root of unity,
we conclude that $\beta=\gamma=0$.
Thus, $L_E(E)$ is a scalar multiple of $E$,
and the same argument applies to $L_F(E)$.
Similarly, by \eqref{eq:skew-1Kinv},
we may write $L_E(F)=\alpha F+\beta EK^{-1}
+\gamma(1-K^{-1}).$
Using $KL_E(F)K^{-1}=q^{-2}L_E(F)$,
we obtain
\[
\begin{aligned}
q^{-2}\alpha F+q^2\beta EK^{-1}
+\gamma(1-K^{-1})=q^{-2}\alpha F+q^{-2}\beta EK^{-1}
+q^{-2}\gamma(1-K^{-1}).
\end{aligned}
\]
Again, linear independence and the assumption
on $q$ imply $\beta=\gamma=0$.
Thus, $L_E(F)$ is a scalar multiple of $F$,
and the same holds for $L_F(F)$.
Consequently,
\[
L_E(E)=uE,\quad L_F(E)=vE,\quad
L_E(F)=wF,\quad L_F(F)=xF
\]
for some $u,v,w,x\in\mathbb{C}$.
Applying \eqref{eq:twisted-leibniz-E}
and \eqref{eq:twisted-leibniz-F} to the
commutator relation in \eqref{eq:uq-rel},
and using $L_E(K^{\pm1})=L_F(K^{\pm1})=0$,
we obtain
\[
\begin{aligned}
0=L_E(EF-FE)
&=(u+tw)(EF-t^{-1}FE),\\
0=L_F(EF-FE)
&=(tv+x)(EF-t^{-1}FE).
\end{aligned}
\]
Since $EF-t^{-1}FE\neq0$, we conclude that
$w=-ut^{-1}$ and $x=-tv$, and thus
\begin{equation}\label{eq:generator-actions}
\begin{aligned}
L_E(E)&=uE, \qquad
L_E(F)&=-ut^{-1}F, \qquad
L_F(E)&=vE, \qquad
L_F(F)&=-tvF.
\end{aligned}
\end{equation}
We now determine all possible triples
$(t,u,v)$ and prove that each of them uniquely
determines the post-Hopf structure.  Applying \eqref{eq:PH2}, together with
the relations $KE=q^2EK$ and $KF=q^{-2}FK$,
we obtain
\[
\begin{aligned}
\Phi_tL_E
&=L_KL_E=tL_{KE}
=tq^2L_{EK}=tq^2L_E\Phi_t,\\
\Phi_tL_F
&=L_KL_F=t^{-1}L_{KF}
=t^{-1}q^{-2}L_{FK}
=t^{-1}q^{-2}L_F\Phi_t.
\end{aligned}
\]
Evaluating these identities at $E$ and using
\eqref{eq:generator-actions}, we find $utE=t^2q^2uE$ and $vtE=q^{-2}vE$.
Since $t\ne0$, either $u\ne0$ or $v\ne0$
implies $t=q^{-2}$.

We next derive a relation between $u$ and $v$.
Applying \eqref{eq:PH2} to the pairs $(E,F)$
and $(F,E)$ gives
\[
\begin{aligned}
L_EL_F&=L_{EF}-ut^{-1}L_{KF},\\
L_FL_E&=t^{-1}L_{FE}+vL_E.
\end{aligned}
\]
Moreover, since $L_KL_F=t^{-1}L_{KF}$,
we have $L_{KF}(E)=t^2vE$.
Evaluating the preceding identities at $E$,
we therefore obtain
\[
\begin{aligned}
L_{EF}(E)
&=L_EL_F(E)+ut^{-1}L_{KF}(E)
=uv(1+t)E,\\
L_{FE}(E)
&=tL_FL_E(E)-tvL_E(E)=0.
\end{aligned}
\]
On the other hand, the commutator relation
in \eqref{eq:uq-rel} gives
\[
L_{EF}(E)-L_{FE}(E)
=\frac{L_K(E)-L_{K^{-1}}(E)}{q-q^{-1}}
=\frac{t-t^{-1}}{q-q^{-1}}E,
\]
and thus
\begin{equation}\label{eq:uv-condition}
uv(1+t)=\frac{t-t^{-1}}{q-q^{-1}}.
\end{equation}

If $u=v=0$, equation \eqref{eq:uv-condition}
gives $t=\pm1$. If at least one of $u,v$
is nonzero, we have $t=q^{-2}$, and
\eqref{eq:uv-condition} yields $uv=-q$.
In particular, both $u$ and $v$ must
be nonzero. Thus, the only possibilities are
\begin{equation}\label{eq:cases}
(t,u,v)=(1,0,0),\qquad
(t,u,v)=(-1,0,0),\qquad
t=q^{-2},\quad uv=-q.
\end{equation}
The triple $(t,u,v)$ uniquely determines the
post-Hopf structure. Namely, we have
$L_{K^{\pm1}}=\Phi_t^{\pm1}$, while
\eqref{eq:twisted-leibniz-E}--\eqref{eq:generator-actions} uniquely determine $L_E$ and
$L_F$ on the whole algebra. By induction on
products of generators, both operators act
by scalars depending only on the degree.
Consequently, applying \eqref{eq:PH2}
successively to homogeneous elements
determines $L_a$ for every PBW basis
element $a=F^rK^kE^s$.
One can check that the operations
$\act_+$, $\act_-$ and $\act_q$
constructed above correspond respectively
to the triples
$(1,0,0)$, $(-1,0,0)$ and
$(q^{-2},1,-q)$.
Moreover, 
any Hopf automorphism $\Phi_z$ sends
the triple $(t,u,v)$ into
$(t,z^{-1}u,zv)$.
In particular, the triples $(1,0,0)$
and $(-1,0,0)$ are fixed by every
$\Phi_z$, and choosing $z=u$
sends the triple $(t,u,v)$ into $(q^{-2},1,-q)$,
so the corresponding post-Hopf structure
is isomorphic to $\act_q$.
Therefore every post-Hopf structure is
isomorphic to one of
$\act_+$, $\act_-$ or $\act_q$. 
Finally, since the parameter $t$ is invariant
under post-Hopf isomorphisms, and $q$ is not a root of unity (hence
$1,-1,q^{-2}$ are pairwise distinct), the three structures are pairwise
non-isomorphic.
\end{proof}

\section{Post-Hopf QUEAs and post-Lie bialgebras}

\subsection{Post-Lie algebras and Lie bialgebras}
\label{subsec:post-lie}

Recall that a \emph{post-Lie algebra} is a Lie algebra
$(\mathfrak g,[\cdot,\cdot])$ equipped with a bilinear
map $\tri:\mathfrak g\otimes\mathfrak g
\to\mathfrak g$ satisfying
\begin{align}
x\tri[y,z]
&=[x\tri y,z]+[y,x\tri z],
\label{eq:postlie-1}\\
[x,y]\tri z
&=x\tri(y\tri z)-y\tri(x\tri z)
-(x\tri y)\tri z+(y\tri x)\tri z
\label{eq:postlie-2}
\end{align}
for all $x,y,z\in\mathfrak g$.
In particular, the map $[x,y]_{\tri}=[x,y]+x\tri y-y\tri x$
defines another Lie bracket on $\mathfrak g$, called the subadjacent Lie bracket, and the map 
\begin{equation}
\label{eq:derivations}
D: (\mathfrak{g},[\cdot,\cdot]_{\tri}) \to \mathrm{Der}(\mathfrak{g},[\cdot,\cdot]), \qquad x \mapsto D_x: \g \to \g, \qquad y \mapsto x \tri y
\end{equation}
is a Lie algebra morphism.
Recall also that a \emph{Rota--Baxter operator of weight one}
on a Lie algebra $(\mathfrak g,[\cdot,\cdot])$ is a
linear map $\varphi:\mathfrak g\to\mathfrak g$ satisfying $[\varphi(x),\varphi(y)]
=\varphi\bigl(
[\varphi(x),y]+[x,\varphi(y)]+[x,y]
\bigr)$
for all $x,y\in\mathfrak g$.

The relationship between post-Lie and post-Hopf
algebras is given by the following result
\cite{LST,RSW}.

\begin{proposition}\label{prop:postlie-posthopf}
Let $(A,\act)$ be a post-Hopf algebra.
Then the space of primitive elements $P(A)$ is a post-Lie algebra with
the commutator bracket and the restriction
of $\act$.
Conversely, every post-Lie algebra
$(\mathfrak g,[\cdot,\cdot],\tri)$ induces a
unique post-Hopf structure on its universal
enveloping algebra $U(\mathfrak g)$
whose restriction to
$\mathfrak g\otimes\mathfrak g$
coincides with $\tri$. These constructions are functorial and
establish an adjunction between post-Lie
algebras and post-Hopf algebras, with the
universal enveloping algebra functor left
adjoint to the primitive-elements functor.
\end{proposition}

Next, recall that a \emph{Lie bialgebra} is a Lie algebra
$(\mathfrak g,[\cdot,\cdot])$ equipped with a linear map
$\delta:\mathfrak g\to\Lambda^2\mathfrak g$, called the
\emph{Lie cobracket}, satisfying the co-Jacobi identity
and the cocycle condition
\begin{equation}\label{eq:lie-bialgebra}
\delta([x,y])
=(\operatorname{ad}_{x}\otimes\id
+\id\otimes\operatorname{ad}_{x})
\delta(y)
-(\operatorname{ad}_y\otimes\id
+\id\otimes\operatorname{ad}_y)
\delta(x)
\end{equation}
for all $x,y\in\mathfrak g$.

In the sequel, we investigate the compatibility
between a post-Lie product and a Lie cobracket
on the same Lie algebra. We derive a natural
compatibility condition from semiclassical
limits of post-Hopf quantized universal
enveloping algebras, leading to the notion
of a \emph{post-Lie bialgebra}.

\subsection{Post-Lie bialgebras as semiclassical limits of post-Hopf QUEAs}
\label{subsec:semiclassical}

From now on, we fix a formal parameter $\hbar$. By \emph{topological Hopf algebra} we mean a Hopf monoid in the category of topologically free modules over $\mathbb{C}[[\hbar]]$, and the completed tensor
product of topologically free $\mathbb{C}[[\hbar]]$-modules will also be denoted by $\otimes$.

Recall that a \emph{quantized universal enveloping
algebra} (QUEA) is a topological Hopf algebra
$A_\hbar$ such that
$A_\hbar/\hbar A_\hbar\cong U(\mathfrak g)$
as Hopf algebras, for some Lie algebra
$\mathfrak g$. The Lie algebra $\mathfrak{g}= P(A_\hbar/\hbar A_\hbar)$ is called the \emph{semiclassical limit} of $A_\hbar$.
It is well known that the semiclassical limit of
a QUEA carries a natural Lie bialgebra structure
\cite{Drinfeld,ChariPressley}. More precisely,
for $x\in\mathfrak g$, choose a lift
$\widetilde{x}\in A_\hbar$ and define
\begin{equation}\label{eq:semiclassical-cobracket}
\delta(x)=
\frac{\Del_\hbar(\widetilde{x})
-\Del_\hbar^{\mathrm{op}}(\widetilde{x})}
{\hbar}  \ \mod \hbar.
\end{equation}
This expression is independent of the chosen lift
and defines a Lie bialgebra cobracket
$\delta:\mathfrak g\to\Lambda^2\mathfrak g$.

\begin{proposition}\label{prop:posthopf-limit}
Let $A_\hbar$ be a QUEA with $A_\hbar/\hbar A_\hbar\cong U(\mathfrak g)$ and suppose that there is a topological post-Hopf structure $\act_\hbar$ on $A_\hbar$. Then the semiclassical limit $\mathfrak g$ of 
$A_\hbar$ carries a natural post-Lie algebra
structure.
\end{proposition}

\begin{proof}
Specializing $\act_\hbar$ at $\hbar=0$
gives a post-Hopf structure $\act_0$ on
$U(\mathfrak g)$.
We first show that $\act_0$ restricts
to $\mathfrak g$. Since $\act_0$ is a coalgebra
morphism, for $x,y\in\mathfrak g$ we have
\[
\begin{aligned}
\Del_0(x\act_0 y)
=(x_1\act_0 y_1)\otimes
  (x_2\act_0 y_2)&=(x\act_0 y)\otimes1
 +(x\act_0 1)\otimes y
 +y\otimes(x\act_0 1)
 +1\otimes(x\act_0 y)\\
 &=(x\act_0 y)\otimes1
 +1\otimes(x\act_0 y),
\end{aligned}
\]
where we used that $x,y \in P(U(\mathfrak{g}))$,
together with \eqref{eq:norm}.
Thus, $x\act_0 y$ is primitive, and $x\tri y:=x\act_0 y$
defines a bilinear map on $\mathfrak g$.
We now verify \eqref{eq:postlie-1},\eqref{eq:postlie-2}.
By \eqref{eq:PH1}
and $\Del_0(x)=x\otimes1+1\otimes x$,
we obtain
\[
\begin{aligned}
x\act_0(ab)
=(x\act_0 a)(1\act_0 b)
 +(1\act_0 a)(x\act_0 b)=(x\act_0 a)b+a(x\act_0 b)
\end{aligned}
\]
for all $x\in\mathfrak g$ and
$a,b\in U(\mathfrak g)$. Consequently, for $y,z\in\mathfrak g$,
\[
\begin{aligned}
x\tri[y,z]
=x\act_0(yz-zy)=(x\tri y)z+y(x\tri z)
 -(x\tri z)y-z(x\tri y)=[x\tri y,z]+[y,x\tri z],
\end{aligned}
\]
which proves \eqref{eq:postlie-1}.
Similarly, by \eqref{eq:PH2}
we obtain
\[
\begin{aligned}
x\tri(y\tri z)
=x\act_0(y\act_0 z)=\bigl(x(1\act_0 y)
       +1(x\act_0 y)\bigr)\act_0 z=(xy+x\tri y)\act_0 z
\end{aligned}
\]
for all $x,y,z \in \mathfrak{g}$, hence
\[
\begin{aligned}
[x,y]\tri z
=(xy-yx)\act_0 z=x\tri(y\tri z)-y\tri(x\tri z) -(x\tri y)\tri z
 +(y\tri x)\tri z,
\end{aligned}
\]
which is \eqref{eq:postlie-2}.
\end{proof}
Combining the last result with the established theory of QUEAs, we obtain that the semiclassical limit of a post-Hopf QUEA
carries both a post-Lie algebra structure
$(\mathfrak g,[\cdot,\cdot],\tri)$ and a Lie
bialgebra structure
$(\mathfrak g,[\cdot,\cdot],\delta)$.
These structures are not independent, as shown
by the following result.

\begin{theorem}\label{thm:semiclassical-compatibility}
Let $A_\hbar$ be a post-Hopf QUEA, and let
$\tri$ and $\delta$ be the induced post-Lie product
and Lie cobracket on $\mathfrak g$. Then for all $x,y\in\mathfrak g$ we have
\begin{equation}\label{eq:postlie-compatibility}
\delta(x\tri y)
=(D_x\otimes\id+\id\otimes D_x)\delta(y).
\end{equation}
\end{theorem}

\begin{proof}
Let $x,y\in\mathfrak g$, and choose lifts
$\widetilde x,\widetilde y\in A_\hbar$. Since $y$ is primitive
in $U(\mathfrak g)$, there is an element $B_y\in A_\hbar\ten A_\hbar$ such that $\Del_\hbar(\widetilde y)
=\widetilde y\ten1+1\ten\widetilde y
+\hbar B_y$.
Denoting by $\overline{B}_y$ its specialization
at $\hbar=0$, we have $\delta(y)=\overline{B}_y-\tau(\overline{B}_y),$
where $\tau$ denotes the flip of tensor factors.
Since $\act_\hbar$ is a coalgebra morphism,
writing $B_y=\sum_i b_i\ten c_i$ gives
\[
\begin{aligned}
\Del_\hbar(\widetilde x\act_\hbar\widetilde y)
=(\widetilde x_1\act_\hbar\widetilde y_1)
  \ten(\widetilde x_2\act_\hbar\widetilde y_2)=(\widetilde x\act_\hbar\widetilde y)\ten1
  +1\ten(\widetilde x\act_\hbar\widetilde y)+\hbar\sum_i
  (\widetilde x_1\act_\hbar b_i)
  \ten(\widetilde x_2\act_\hbar c_i),
\end{aligned}
\]
where the second equality follows from \eqref{eq:norm} and from the fact that $\varepsilon$ is a counit.
Since $\widetilde x$ is primitive modulo $\hbar$,
specializing the last term at $\hbar=0$ yields
\[
\sum_i
\bigl((x\act_0\overline b_i)\otimes\overline c_i
+\overline b_i\otimes
 (x\act_0\overline c_i)\bigr)
=(D_x\otimes\id+\id\otimes D_x)\overline B_y,
\]
where $D_x(a)=x\act_0 a$ is the derivation
of $U(\mathfrak g)$ extending
$D_x(y)=x\tri y$ on $\mathfrak g$,
as shown in the proof of
Proposition~\ref{prop:posthopf-limit}.
Therefore, we have
\[
\left.
\frac{
\Del_\hbar(\widetilde x\act_\hbar\widetilde y)
-(\widetilde x\act_\hbar\widetilde y)\ten1
-1\ten(\widetilde x\act_\hbar\widetilde y)}
{\hbar}
\right|_{\hbar=0}
=(D_x\otimes\id+\id\otimes D_x)\overline B_y.
\]
Since $\widetilde x\act_\hbar\widetilde y$
is a lift of $x\tri y$, applying
$\id-\tau$ to both sides yields
\[
\begin{aligned}
\delta(x\tri y)
=(\id-\tau)
  (D_x\otimes\id+\id\otimes D_x)
  \overline B_y
=(D_x\otimes\id+\id\otimes D_x)\delta(y).
\end{aligned}
\]
\end{proof}

Theorem~\ref{thm:semiclassical-compatibility}
motivates the following definition.

\begin{definition}\label{def:postlie-bialgebra}
A \emph{post-Lie bialgebra} is a quadruple
$(\mathfrak g,[\cdot,\cdot],\tri,\delta)$
such that $(\mathfrak g,[\cdot,\cdot],\tri)$
is a post-Lie algebra,
$(\mathfrak g,[\cdot,\cdot],\delta)$
is a Lie bialgebra, and
\eqref{eq:postlie-compatibility} holds
for all $x,y\in\mathfrak g$.
\end{definition}

By Theorem~\ref{thm:semiclassical-compatibility},
the semiclassical limit of every post-Hopf QUEA
carries a natural post-Lie bialgebra structure.
This construction extends the usual semiclassical
limit functor (from QUEAs to Lie bialgebras) to post-Hopf
QUEAs and post-Lie bialgebras.
Note also that the compatibility condition \eqref{eq:postlie-compatibility}  has the following natural
interpretation: the map \eqref{eq:derivations} corestricts to a Lie algebra morphism
\[
D:(\mathfrak g,[\cdot,\cdot]_{\tri})
\longrightarrow
\operatorname{BiDer}(\mathfrak g,[\cdot,\cdot],\delta),
\]
where
$\operatorname{BiDer}(\mathfrak g,[\cdot,\cdot],\delta)$
denotes the Lie algebra of derivations of the
Lie bracket that are also coderivations of the
Lie cobracket.

\begin{remark}
Our notion of post-Lie bialgebra arises from
semiclassical limits of post-Hopf QUEAs.
It differs from the pp-post-Lie bialgebras
introduced by Lu--Bai--Guo \cite{LuBaiGuo},
which are defined through Manin triples
and partially split post-Lie structures. Note also that our notion of post-Lie bialgebra is not self-dual.
\end{remark}

\subsection{A post-Hopf structure on
$U_\hbar(\mathfrak{sl}_2)$ and post-Lie bialgebra structures on $\mathfrak{sl}_2$}
\label{subsec:hbar}

We now construct a nontrivial topological
post-Hopf structure on $U_\hbar(\mathfrak{sl}_2)$
and determine the corresponding post-Lie
bialgebra.
Recall that $U_\hbar(\mathfrak{sl}_2)$
is the
topologically free $\mathbb{C}[[\hbar]]$-algebra
generated by $E,F,H$, subject to
\begin{equation}\label{eq:hbar-rel}
[H,E]=2E,\qquad
[H,F]=-2F,\qquad
[E,F]=
\frac{e^{\hbar H}-e^{-\hbar H}}
{e^\hbar-e^{-\hbar}}.
\end{equation}

Its Hopf algebra structure is given by comultiplication
\begin{equation}\label{eq:hbar-cop}
\begin{aligned}
\Del_\hbar(H)
=H\ten1+1\ten H,\quad
\Del_\hbar(E)
=E\ten1+e^{\hbar H}\ten E,\quad
\Del_\hbar(F)
=F\ten e^{-\hbar H}+1\ten F,
\end{aligned}
\end{equation} 
counit
$\eps_\hbar(E)=\eps_\hbar(F)
=\eps_\hbar(H)=0$, and antipode $S_\hbar(H)=-H$, 
$S_\hbar(E)=-e^{-\hbar H}E$,  and
$S_\hbar(F)=-Fe^{\hbar H}$.
Clearly,  $U_\hbar(\mathfrak{sl}_2)/
\hbar U_\hbar(\mathfrak{sl}_2)
\cong U(\mathfrak{sl}_2)$
as Hopf algebras. 
Moreover, the degree assignment introduced in
\eqref{eq:grading} is compatible with the
$\hbar$-adic completion, with
$\deg(E)=1$, $\deg(F)=-1$, and $\deg(H)=0$.

\begin{theorem}\label{thm:formal}
There exists a nontrivial topological post-Hopf
structure $\act_\hbar$ on
$U_\hbar(\mathfrak{sl}_2)$ whose values on
generators are
\begin{equation}\label{eq:formal-gen}
\begin{array}{c|ccc}
\act_\hbar & E & F & H\\ \hline
E & E & -e^{2\hbar}F & 0\\
F & -e^\hbar E & e^{-\hbar}F & 0\\
H & -2E & 2F & 0
\end{array}
\end{equation}
\end{theorem}

\begin{proof}
The proof consists in adapting the construction of Subsection~1.3.2
to the formal setting. Consider the
$\hbar$-adically complete algebra
$V_\hbar=\mathbb{C}[z,z^{-1}][[\hbar]]$,
with $v_j=z^j$. The formulas
\[
\begin{aligned}
H\cdot v_j&=2jv_j, \qquad
E\cdot v_j
=\frac{e^{2j\hbar}-1}{e^{-2\hbar}-1}\,v_{j+1},\qquad
F\cdot v_j
=\frac{1-e^{-2j\hbar}}{e^\hbar-e^{-\hbar}}\,v_{j-1}
\end{aligned}
\]
define a $U_\hbar(\mathfrak{sl}_2)$-module
algebra structure on $V_\hbar$.
As in Lemma~\ref{lem:chi}, this action
determines continuous functionals
$\chi_m:U_\hbar(\mathfrak{sl}_2)
\to\mathbb{C}[[\hbar]]$, defined as in \eqref{eq:chi-def} and satisfying \eqref{eq:chi-axioms}.
These identities remain valid with completed
tensor products, so the construction of
Proposition~\ref{prop:graded} extends to a topological post-Hopf structure
defined, for homogeneous
$b\in U_\hbar(\mathfrak{sl}_2)_n$, by $a\act_\hbar b=\chi_n(a)b.$
Evaluating the functionals on the generators
gives \eqref{eq:formal-gen}. 
\end{proof}

We now determine the semiclassical limit of the
post-Hopf structure constructed in
Theorem~\ref{thm:formal}.
Let $e,f,h$ denote the images of $E,F,H$,
respectively, modulo $\hbar$.
The induced Lie cobracket on $\mathfrak{sl}_2$
is the standard one, given by
\begin{equation}\label{eq:sl2-cobracket}
\delta(h)=0,\qquad
\delta(e)=h\wedge e,\qquad
\delta(f)=h\wedge f.
\end{equation}
On the other hand, reducing
\eqref{eq:formal-gen} modulo $\hbar$
shows that the induced post-Lie product
is given by $x\tri y=[\varphi(x),y]$, where $\varphi$ is the Rota-Baxter operator of weight one determined by
\begin{equation}\label{eq:phi}
\varphi(e)=\frac12h,\qquad
\varphi(f)=-\frac12h,\qquad
\varphi(h)=-h.
\end{equation}
Thus, by
Theorem~\ref{thm:semiclassical-compatibility},
the quadruple
$(\mathfrak{sl}_2,[\cdot,\cdot],\tri,\delta)$
is a post-Lie bialgebra.
\begin{remark}
The post-Lie algebra structure \eqref{eq:phi}
belongs to the 
isomorphism class of the one associated to the Rota-Baxter operator of weight one determined by
\[
\varphi_*(e)=0,\qquad
\varphi_*(f)=-f-h,\qquad
\varphi_*(h)=0,
\]
see \cite[Subsection~1.1.2]{RSW}.
A direct computation gives $\varphi=\Psi \circ \varphi_* \circ\Psi^{-1},$ where $\Psi = \mathrm{Ad}(G)$, and
\[
G=
\begin{pmatrix}
1&0\\
-1&2
\end{pmatrix}.
\]
\end{remark}

Next, recall that post-Lie structures on complete Lie
algebras are in one-to-one correspondence with
Rota--Baxter operators of weight one
(see \cite[Remark~3]{RSW}).
For simple Lie algebras, this correspondence
extends to post-Lie bialgebras as follows.

\begin{lemma}\label{lem:rb}
Let $(\mathfrak g,[\cdot,\cdot],\delta)$ be a
finite-dimensional complex simple Lie bialgebra.
Then compatible post-Lie products on
$\mathfrak g$ are in one-to-one correspondence
with weight-one Rota--Baxter operators
$\varphi:\mathfrak g\to\mathfrak g$
satisfying $\delta\circ\varphi=0$.
\end{lemma}

\begin{proof}
Since $\mathfrak g$ is simple, it is complete.
Hence every post-Lie product is uniquely of the form
$x\tri y=[\varphi(x),y]$, where $\varphi$ is a
weight-one Rota--Baxter operator.
By the Lie bialgebra cocycle identity
\eqref{eq:lie-bialgebra}, we have
\[
\begin{aligned}
\delta(x\tri y)
&=(\operatorname{ad}_{\varphi(x)}\otimes\id
+\id\otimes\operatorname{ad}_{\varphi(x)})
\delta(y)-
(\operatorname{ad}_y\otimes\id
+\id\otimes\operatorname{ad}_y)
\delta(\varphi(x)).
\end{aligned}
\]
Thus, \eqref{eq:postlie-compatibility} holds
if and only if
$\delta(\varphi(x))\in
(\Lambda^2\mathfrak g)^{\mathfrak g}$
for every $x\in\mathfrak g$.
Since $\mathfrak g$ is simple,
$(\Lambda^2\mathfrak g)^{\mathfrak g}=0$,
and the claim follows.
\end{proof}

We now apply Lemma~\ref{lem:rb} to classify
the post-Lie structures compatible with the
standard Lie cobracket on $\mathfrak{sl}_2$.

\begin{theorem}\label{thm:sl2-plb}
For the standard Lie bialgebra
$(\mathfrak{sl}_2,[\cdot,\cdot],\delta)$
of \eqref{eq:sl2-cobracket}, every compatible
post-Lie product is induced either by
$\varphi=0$ or by
\begin{equation}\label{eq:phi-family}
\varphi_a(e)=ah,\qquad
\varphi_a(f)=-\frac{1}{4a}h,\qquad
\varphi_a(h)=-h,
\qquad a\in\mathbb{C}^\times,
\end{equation}
through $x\tri y=[\varphi(x),y]$.
All products in \eqref{eq:phi-family}
are isomorphic by Lie bialgebra automorphisms.
Consequently, there are exactly two
isomorphism classes of compatible post-Lie
products.
\end{theorem}

\begin{proof}
By Lemma~\ref{lem:rb}, compatible post-Lie
products correspond to weight-one
Rota--Baxter operators $\varphi$ satisfying
$\delta\circ\varphi=0$.
Since $\ker\delta=\mathbb{C} h$, we may write
\[
\varphi(e)=ah,\qquad
\varphi(f)=bh,\qquad
\varphi(h)=ch.
\]
The Rota--Baxter identity applied to the pairs
$(e,f)$, $(h,e)$ and $(h,f)$ yields
\[
c-4ab=0,\qquad
a(1+c)=0,\qquad
b(1+c)=0.
\]
Thus, either $a=b=c=0$, or
$c=-1$ and $ab=-1/4$, giving
\eqref{eq:phi-family}.
For every $t\in\mathbb{C}^\times$, the map
\[
\Psi_t(e)=te,\qquad
\Psi_t(f)=t^{-1}f,\qquad
\Psi_t(h)=h
\]
is a Lie bialgebra automorphism of
$(\mathfrak{sl}_2,[\cdot,\cdot],\delta)$.
Since
$\Psi_t \circ \varphi_a \circ \Psi_t^{-1}=\varphi_{a/t}$,
all the nonzero products are isomorphic.
The zero product forms a separate
isomorphism class, completing the proof.
\end{proof}

Note that the post-Lie product obtained as the semiclassical
limit of Theorem~\ref{thm:formal} corresponds
to $\varphi_{1/2}$.

\subsection{The quantization problem for post-Lie bialgebras}
\label{subsec:quantization}

The semiclassical functor arising from
Theorem~\ref{thm:semiclassical-compatibility}
naturally leads to the converse problem
of quantizing post-Lie bialgebras.

\begin{definition}
A \emph{quantization} of a post-Lie bialgebra $(\mathfrak g,[\cdot,\cdot],\tri,\delta)$
is a post-Hopf QUEA $(A_\hbar,\act_\hbar)$,
together with a Hopf algebra isomorphism $A_\hbar/\hbar A_\hbar\cong U(\mathfrak g),$
such that, under this isomorphism, the induced Lie cobracket is
$\delta$ and the post-Lie product induced
by the specialization at zero of $\act_\hbar$ coincides
with $\tri$.
\end{definition}

The \emph{quantization problem for post-Lie
bialgebras} asks whether every post-Lie
bialgebra admits such a quantization.
It is well known that every Lie bialgebra admits
a QUEA quantization
\cite{EtingofKazhdan,Severa,RivezziSchnitzer}.
However, it is not yet clear how to extend
these quantization constructions to
arbitrary post-Lie bialgebras.
In the case of $\mathfrak{sl}_2$,
Theorems~\ref{thm:formal} and
\ref{thm:sl2-plb} show that both isomorphism
classes of post-Lie bialgebra structures
with the standard Lie bracket and cobracket
admit quantizations, the trivial class
being realized by the trivial post-Hopf
structure. The existence of quantizations
for arbitrary post-Lie bialgebras remains
an open problem.

\textbf{Contacts:}\\
~\\
arivezzi.math@pmf.hr\\
\url{https://sites.google.com/view/andrearivezzipersonalwebpage/}

\end{document}